\documentclass[11pt]{article}
\usepackage{}
\usepackage{mathrsfs}
\usepackage{amsfonts}

\usepackage{CJK}
\usepackage{amsfonts,amssymb,amsmath,mathrsfs}

\usepackage{color,latexsym,amsfonts}
 \newtheorem{theorem}{Theorem}[section]
 \newtheorem{lemma}[theorem]{Lemma}
 \newtheorem{corollary}[theorem]{Corollary}
 \newtheorem{proposition}[theorem]{Proposition}
 
 \newtheorem{Definition}[theorem]{Definition}
 \newtheorem{remark}[theorem]{Remark}
 \newtheorem{condition}[theorem]{Condition}

 \def\blemma{\begin{lemma}\sl{}\def\elemma{\end{lemma}}}
 \def\btheorem{\begin{theorem}\sl{}\def\etheorem{\end{theorem}}}
 \def\bcorollary{\begin{corollary}\sl{}\def\ecorollary{\end{corollary}}}

 \def\bremark{\begin{remark}\sl{}\def\eremark{\end{remark}}}

 \def\beqlb{\begin{eqnarray}}\def\eeqlb{\end{eqnarray}}
 \def\beqnn{\begin{eqnarray*}}\def\eeqnn{\end{eqnarray*}}

 \def\<{\langle}\def\>{\rangle}

 \def\eqref#1{{\rm(\ref{#1})}}

\newcommand{\ra}{\rightarrow}

\def\D{\textup{D}}

\def\d{\textup{d}}
\def\I{\textup{I}}
\def\e{\textup{e}}

\def\fin{\hfill$\square$}

\def\newdot{{\kern.8pt\cdot\kern.8pt}}

\def\R{\mathbb{R}}
\def\E{\mathbb{E}}
\def\P{\mathbb{P}}
\def\D{\mathbb{D}}

\def\<{\langle}
\def\>{\rangle}
\def\Proof.{\noindent{\bf Proof.}}

\begin{document}

\

\noindent{}

\bigskip\bigskip

\centerline{\Large\bf Bismut formulae and applications for}

\smallskip

\centerline{\Large\bf
stochastic (functional) differential equations driven by}

\smallskip

\centerline{\Large\bf
fractional Brownian motions\footnote{Supported by the Research project of Natural Science Foundation of
   Anhui Provincial Universities (Grant No. KJ2013A134), National Natural Science Foundation of China
(Grant No. 11371029).}}

\smallskip
\
\bigskip\bigskip

\centerline{Xiliang Fan}

\bigskip

\centerline{Department of Statistics, Anhui Normal University,}

\centerline{Wuhu 241003, China}

\smallskip

\bigskip\bigskip

{\narrower{\narrower

\noindent{\bf Abstract.} By using Malliavin calculus, Bismut derivative formulae are established
for a class of stochastic (functional) differential equations driven by fractional Brownian motions.
As applications, Harnack type inequalities and strong Feller property are presented.
}}


\bigskip



\section{Introduction}

\setcounter{equation}{0}

The Bismut formula, initiated in  \cite{Bismut84},
has become a very effective tool in stochastic analysis.
Since then, it has been well studied in the setting of Markov process.
For instance, the readers may refer to \cite{Dong&Xie10,Elworthy&Li94} (using martingale method),
\cite{Bao&Wang&Yuan13b,Guillin&Wang12,Wang&Xu10a} (via coupling argument),
\cite{Bao&Wang&Yuan13c,Priola06,Wang&Zhang13a,Zhang10a} (by Malliavin calculus),
and references therein.

In this article, we are interested in stochastic differential (functional) equations driven by fractional Brownian motions.
The results on the existence and uniqueness of solutions of such stochastic equations were obtained using mainly the theory of rough path analysis introduced in \cite{Lyons98a} or a fractional integration by parts formula (see \cite{Zahle98a}).
The readers may refer to \cite{Boufoussi&Hajji11,Coutin&Qian02a,Diopa&Garrido-Atienza14,Ferrante&Rovira06,Neuenkirch&Nourdin&Tindel08,Nualart&Ouknine02b,Nualart&Rascanu02a} and references therein.
In \cite{Baudoin&Ouyang11,Nourdin&Simon06,Nualart&Saussereau09}, the authors studied the regularities of the solutions.
Hairer and Pillai \cite{Hairer05,Hairer&Pillai11} investigated the ergodicity of the solution,
and the convergence rate toward the stationary solution.
Saussereau \cite{Saussereau12} proved Talagrand's transportation inequalities for the law of the solution and studied the asymptotic behaviors.
In \cite{Baudoin&Ouyang&Tindel11a}, the authors showed the logarithmic Sobolev inequalities for the law of the solution at a fixed time.
This work is strongly motivated by the study of strong Feller property for operators associated with the solutions in a non-Markovian context.
To the best of our knowledge, little seems to be known on this subject.
Recently, by using the coupling argument the author \cite{Fan13,Fan14} established Harnack type inequalities (which is stronger than strong Feller property) for SDEs (without delay) with $H<1/2$ and $H>1/2$, respectively.
Inspired by the work \cite{Bao&Wang&Yuan13c}, where derivative formulae are given for functional SPDEs driven by Brownian motion,
we will show the strong Feller property and Harnack type inequalities as consequences of stronger property: Bismut formulae established by Malliavin calculus, for SDEs and SFDEs driven by fractional Brownian motions.

The rest of this paper is organized as follows.
In the next section we present some preliminaries on fractional calculus and fractional Brownian motion.
Section 3 is devoted to the case of SDE.
Then, in Section 4 we consider the SFDE case.

\section{Preliminaries}

\setcounter{equation}{0}

\subsection{Fractional calculus}
In the part, we recall some basic facts about fractional calculus.
An exhaustive survey can be found in \cite{Samko&Kilbas&Marichev}.

Let $a,b\in\R$ and $a<b$.
For $f\in L^1([a,b],\R)$ and $\alpha>0$, the left-sided (resp. right-sided) fractional Riemann-Liouville integral of $f$ of order $\alpha$
on $[a,b]$ is given by
\beqnn
I_{a+}^\alpha f(x)=\frac{1}{\Gamma(\alpha)}\int_a^x\frac{f(y)}{(x-y)^{1-\alpha}}\d y\\
\left(\mbox{resp.}\ \ I_{b-}^\alpha f(x)=\frac{(-1)^{-\alpha}}{\Gamma(\alpha)}\int_x^b\frac{f(y)}{(y-x)^{1-\alpha}}\d y\right),
\eeqnn
where $x\in(a,b)$ a.e., $(-1)^{-\alpha}=\e^{-i\alpha\pi},\Gamma$ denotes the Euler function.\\
They extend the usual $n$-order iterated integrals of $f$ for $\alpha=n\in\mathbb{N}$.
By the definitions, we can get the first composition formulae
\beqnn
I_{a+}^\alpha(I_{a+}^\beta f)=I_{a+}^{\alpha+\beta}f,\quad I_{b-}^\alpha(I_{b-}^\beta f)=I_{b-}^{\alpha+\beta}f.
\eeqnn

Fractional differentiation may be introduced as an inverse operation.
Let $\alpha\in(0,1)$ and $p\geq1$.
If $f\in I_{a+}^\alpha(L^p([a,b],\R))$ (resp. $I_{b-}^\alpha(L^p([a,b],\R)))$, the function $\phi$ satisfying $f=I_{a+}^\alpha\phi$ (resp. $f=I_{b-}^\alpha\phi$) is unique in $L^p([a,b],\R)$ and it coincides with the left-sided (resp. right-sided) Riemann-Liouville derivative
of $f$ of order $\alpha$ given by
\beqnn
D_{a+}^\alpha f(x)=\frac{1}{\Gamma(1-\alpha)}\frac{\d}{\d x}\int_a^x\frac{f(y)}{(x-y)^\alpha}\d y\quad
\left(\mbox{resp.}\ D_{b-}^\alpha f(x)=\frac{(-1)^{1+\alpha}}{\Gamma(1-\alpha)}\frac{\d}{\d x}\int_x^b\frac{f(y)}{(y-x)^\alpha}\d y\right).
\eeqnn
The corresponding Weyl representation reads as follow
\beqnn
D_{a+}^\alpha f(x)=\frac{1}{\Gamma(1-\alpha)}\left(\frac{f(x)}{(x-a)^\alpha}+\alpha\int_a^x\frac{f(x)-f(y)}{(x-y)^{\alpha+1}}\d y\right)\\
\left(\mbox{resp.}\ \ D_{b-}^\alpha f(x)=\frac{(-1)^\alpha}{\Gamma(1-\alpha)}\left(\frac{f(x)}{(b-x)^\alpha}+\alpha\int_x^b\frac{f(x)-f(y)}{(y-x)^{\alpha+1}}\d y\right)\right),
\eeqnn
where the convergence of the integrals at the singularity $y=x$ holds pointwise for almost all $x$ if $p=1$ and in the $L^p$ sense if $p>1$.

By construction, we have
\beqnn
I_{a+}^\alpha(D_{a+}^\alpha f)=f,\quad \forall f\in  I_{a+}^\alpha(L^p([a,b],\R));\quad D_{a+}^\alpha(I_{a+}^\alpha f)=f,\quad \forall f\in L^1([a,b],\R),
\eeqnn
and moreover there hold the second composition formulae
\beqnn
D_{a+}^\alpha(D_{a+}^\beta f)=D_{a+}^{\alpha+\beta}f,~ f\in I_{a+}^{\alpha+\beta}(L^1([a,b],\R));\quad D_{b-}^\alpha(D_{b-}^\beta f)=D_{b-}^{\alpha+\beta}f, ~f\in I_{b-}^{\alpha+\beta}(L^1([a,b],\R)).
\eeqnn

\subsection{Fractional Brownian motion}
Let $B^H=\{B_t^H, t\in[0,T]\}$ be a $d$-dimensional fractional Brownian motion with Hurst parameter $H\in(0,1)$ defined on the probability
space $(\Omega,\mathscr{F},\mathbb{P})$.
That is, $B^H$ is a centered Gauss process with the covariance function $\E B_t^{H,i}B_s^{H,j}=R_H(t,s)\delta_{i,j}$, where
\beqnn
 R_H(t,s)=\frac{1}{2}\left(t^{2H}+s^{2H}-|t-s|^{2H}\right).
\eeqnn
In particular, if $H=1/2, B^H$ is a $d$-dimensional Brownian motion.
By the above covariance function, one can show that $\E|B_t^{H,i}-B_s^{H,i}|^p=C(p)|t-s|^{pH},\ \forall p\geq 1$.
As a consequence, $B^{H,i}$ have $(H-\epsilon)$-order H\"{o}lder continuous paths for all $\epsilon>0,\ i=1,\cdot\cdot\cdot,d$.

For each $t\in[0,T]$, we denote by $\mathcal {F}_t$ the $\sigma$-algebra generated by the random variables $\{B_s^H:s\in[0,t]\}$ and the
$\mathbb{P}$-null sets.

We denote by $\mathscr{E}$ the set of step functions on $[0,T]$. Let
$\mathcal {H}$ be the Hilbert space defined as the closure of
$\mathscr{E}$ with respect to the scalar product
\beqnn
\langle (I_{[0,t_1]},\cdot\cdot\cdot,I_{[0,t_d]}),(I_{[0,s_1]},\cdot\cdot\cdot,I_{[0,s_d]})\rangle_\mathcal {H}=\sum\limits_{i=1}^dR_H(t_i,s_i).
\eeqnn
The mapping $(I_{[0,t_1]},\cdot\cdot\cdot,I_{[0,t_d]})\mapsto\sum_{i=1}^dB_{t_i}^{H,i}$ can be extended to an isometry between $\mathcal {H}$ and the Gauss space $\mathcal {H}_1$ associated with $B^H$. Denote this isometry by $\phi\mapsto B^H(\phi)$.
On the other hand, from \cite{Decreusefond&Ustunel98a}, we know the
covariance kernel $R_H(t,s)$ can be written as
\beqnn
 R_H(t,s)=\int_0^{t\wedge s}K_H(t,r)K_H(s,r)\d r,
\eeqnn
where $K_H$ is a square integrable kernel given by
\beqnn
K_H(t,s)=\Gamma\left(H+\frac{1}{2}\right)^{-1}(t-s)^{H-\frac{1}{2}}F\left(H-\frac{1}{2},\frac{1}{2}-H,H+\frac{1}{2},1-\frac{t}{s}\right),
\eeqnn
in which $F(\cdot,\cdot,\cdot,\cdot)$ is the Gauss hypergeometric function (for details see \cite{Decreusefond&Ustunel98a} or \cite{Nikiforov&Uvarov88}).\\

Define the linear operator $K_H^*:\mathscr{E}\rightarrow L^2([0,T],\R^d)$ as follows
\beqnn
(K_H^*\phi)(s)=K_H(T,s)\phi(s)+\int_s^T(\phi(r)-\phi(s))\frac{\partial K_H}{\partial r}(r,s)\d r.
\eeqnn
By integration by parts, it is easy to see that when $H>1/2$, the above relation can be rewritten as
\beqnn
(K_H^*\phi)(s)=\int_s^T\phi(r)\frac{\partial K_H}{\partial r}(r,s)\d r.
\eeqnn

By \cite{Alos&Mazet&Nualart01a}, we know that, for all $\phi,\psi\in\mathscr{E}$, $ \langle
K_H^*\phi,K_H^*\psi\rangle_{L^2([0,T],\R^d)}=\langle\phi,\psi\rangle_\mathcal {H}$ holds.
From bounded linear transform theorem, $K_H^*$ can be extended to an isometry between $\mathcal{H}$ and $L^2([0,T],\R^d)$.
Therefore, according to \cite{Alos&Mazet&Nualart01a}, the process $\{W_t=B^H((K_H^*)^{-1}{\rm I}_{[0,t]}),t\in[0,T]\}$ is a
Wiener process, and $B^H$ has the following integral representation
\beqnn
 B^H_t=\int_0^tK_H(t,s)\d W_s.
\eeqnn

According to \cite{Decreusefond&Ustunel98a}, the operator $K_H: L^2([0,T],\mathbb{R}^d)\rightarrow I_{0+}^{H+1/2 }(L^2([0,T],\mathbb{R}^d))$ associated with the kernel $K_H(\cdot,\cdot)$ is defined as follows
\beqnn
 (K_Hf^i)(t)=\int_0^tK_H(t,s)f^i(s)\d s,\ i=1,\cdot\cdot\cdot,d.
\eeqnn
It is an isomorphism and for each $f\in L^2([0,T],\mathbb{R}^d)$,
\beqnn
 (K_H f)(s)=I_{0+}^{2H}s^{1/2-H}I_{0+}^{1/2-H}s^{H-1/2}f,\ H\leq1/2,\\
 (K_H f)(s)=I_{0+}^{1}s^{H-1/2}I_{0+}^{H-1/2}s^{1/2-H}f,\ H\geq1/2.
\eeqnn
 As a consequence, for every $h\in I_{0+}^{H+1/2}(L^2([0,T],\R^d))$, the inverse operator $K_H^{-1}$
is of the following form
\beqlb\label{2.1}
(K_H^{-1}h)(s)=s^{H-1/2}D_{0+}^{H-1/2}s^{1/2-H}h',\ H>1/2,
\eeqlb
\beqlb\label{2.2}
(K_H^{-1}h)(s)=s^{1/2-H}D_{0+}^{1/2-H}s^{H-1/2}D_{0+}^{2H}h,\ H<1/2.
\eeqlb
In particular, if $h$ is absolutely continuous, we get
\beqlb\label{2.3}
 (K_H^{-1}h)(s)=s^{H-1/2}I_{0+}^{1/2-H}s^{1/2-H}h',\ H<1/2.
\eeqlb

The remaining part will be devoted to the Malliavin calculus of fractional Brownian motion.

Let $\Omega$ be the canonical probability space $C_0([0,T],\R^d)$, the set of continuous functions,
null at time $0$, equipped with the supremum norm.
Let $\P$ be the unique probability measure on $\Omega$
such that the canonical process $\{B^H_t; t\in[0,T]\}$ is a $d$-dimensional fractional Brownian motion with Hurst parameter $H$.
Then, the injection $R_H=K_H\circ K_H^*:\mathcal{H}\rightarrow\Omega$ embeds $\mathcal{H}$ densely into $\Omega$ and
$(\Omega,\mathcal{H},\P)$ is an abstract Wiener space in the sense of Gross.
In the sequel we will make this assumption on the underlying probability space.

Let $\mathcal {S}$ denote the set of smooth and cylindrical random variables of the form
\beqnn
F=f(B^H(\phi_1),\cdot\cdot\cdot,B^H(\phi_n)),
\eeqnn
where $n\geq 1, f\in C_b^\infty(\mathbb{R}^n)$, the set of $f$ and all its partial derivatives are bounded, $\phi_i\in\mathcal{H}, 1\leq i\leq n$.
The Malliavin derivative of $F$, denoted by $\mathbb{D}F$, is defined as the $\mathcal {H}$-valued random variable
\beqnn
\mathbb{D}F=\sum_{i=1}^n\frac{\partial f}{\partial x_i}(B^H(\phi_1),\cdot\cdot\cdot,B^H(\phi_n))\phi_i.
\eeqnn
For any $p\geq 1$, we define the Sobolev space $\mathbb{D}^{1,p}$ as the completion of $\mathcal {S}$ with respect to the norm
\beqnn
\|F\|_{1,p}^p=\mathbb{E}|F|^p+\mathbb{E}\|\mathbb{D}F\|^p_{\mathcal {H}}.
\eeqnn
$\delta_H$ is denoted by the divergence operator of $\mathbb{D}$.

\section{Main results for SDEs}

\setcounter{equation}{0}

Consider the following stochastic differential equation (SDE for short) driven by fractional Brownian motion:
\beqlb\label{3.1}
\d X(t)=b(X(t))\d t+\sigma(t)\d B^H(t), \ X(0)=x\in\R^d,
\eeqlb
where $b:\mathbb{R}^d\rightarrow\mathbb{R}^d, \sigma:[0,T]\rightarrow\R^d\times\R^d$ and $H>1/2$.

In this part, we aim to study Bismut formulae for the associated family of Markov operators $(P(t))_{0\leq t\leq T}$:
$$P(t)f(x):=\mathbb{E}f(X^x(t)), \ t\in[0,T], \ f\in \mathscr{B}_b(\mathbb{R}^d),$$
where $X^x(t)$ is the solution to \eqref{3.1} with $X(0)=x$ and $\mathscr{B}_b(\mathbb{R}^d)$ denotes the set of all bounded measurable functions on $\mathbb{R}^d$.
Besides,
for $0<\alpha\leq1$, let $C^\lambda(0,T;\mathbb{R}^d)$ be the space of $\alpha$-H\"{o}lder continuous functions $f:[0,T]\rightarrow\mathbb{R}^d$ and set
$$\|f\|_\alpha:=\sup\limits_{0\leq s<t\leq T}\frac{|f(t)-f(s)|}{|t-s|^\alpha}.$$

To the end, we introduce the following hypothesis: \\
(H1) there exist two positive constants $K$ and $\tilde{K}$ such that
\begin{itemize}
\item[(i)] $|\nabla b(x)|\leq K, \ \forall x\in\R^d$;

\par

\item[(ii)] $\sigma$ is bounded and $|\sigma^{-1}(t)-\sigma^{-1}(s)|\leq\tilde{K}|t-s|^{\alpha_0},\ \forall t,s\in[0,T]$,
where $H-1/2<\alpha_0\leq1$.
\end{itemize}

Below we will give a lemma which is important for the proof of our main results.
\blemma\label{L3.1}
Assume (H1).
Then, there hold $X^i(t)\in\D^{1,2}$ and
\beqnn
\mathbb{D}X^i(t)=\sum_{j=1}^d\int_0^t(\nabla b(X(s)))_{ij}\mathbb{D}X_s^j\d s+\sum\limits_{j=1}^d(\sigma_{ij}\I_{[0,t]})e^j,\  1\leq i\leq d, \ 0\leq t\leq T,
\eeqnn
where $\{e^i\}_{i=1}^d$ is the canonical ONB on $\R^d$.
\elemma
\emph{Proof.}
Obviously, by (H1), it follows from \cite[Theorem 2.1]{Nualart&Rascanu02a} that \eqref{3.1} has a unique solution.
Next, we consider the Picard iteration scheme as follows
\begin{equation}\label{1.1}\nonumber
\left\{
\begin{array}{ll}
\d X^{(n+1)}(t)=b(X^{(n)}(t))\d t+\sigma(t)\d B^H(t), \ n\geq0,\\
X^{(0)}(t)=x.
\end{array} \right.
\end{equation}
Observe that
$$|X^{(n+1)}(t)-X(t)|\leq K\int_0^t|X^{(n)}(s)-X(s)|\d s,$$
then the induction argument implies that
$$|X^{(n+1)}(t)-X(t)|\leq \frac{K^{n+1}}{(n+1)!}(|x|+\|X\|_\infty)t^{n+1}.$$
By  \cite[Theorem 2.1 II]{Nualart&Rascanu02a}, we know that, for any $p\geq1$, there holds $\E \|X\|_\infty^p\<\infty.$
As a consequence, we derive that, for all $p\geq1$,
\beqnn
\E\sup\limits_{0\leq t\leq T}|X^{(n+1)}(t)-X(t)|^p\leq\left[\frac{(T K)^{n+1}}{(n+1)!}\right]^p\cdot\E(|x|+\|X\|_\infty)^p\ra0,\ n\ra\infty.
\eeqnn
On the other hand, we easily get that
\beqnn
\mathbb{D}X^{(n+1),i}(t)=\sum_{j=1}^d\int_0^t(\nabla b(X^{(n)}(s)))_{ij}\mathbb{D}X^{(n),j}(s)\d s+\sum\limits_{j=1}^d(\sigma_{ij}\I_{[0,t]})e^j.
\eeqnn
Again by the induction argument, it follows that
$$\|\mathbb{D}X^{(n+1),i}(t)\|_\mathcal {H}\leq1+Kt+\frac{(Kt)^2}{2}+\cdots+\frac{(Kt)^n}{n!}.$$
Consequently, $\sup_{n\geq0}\|\mathbb{D}X^{(n+1),i}(t)\|_\mathcal {H}\leq\e^{KT}$ holds.
Hence, with the help of \cite[Lemma 1.2.3]{Nualart06a}, we conclude that $X^i(t)\in\D^{1,2}$.
The second assertion follows easily from \eqref{3.1}.
\fin

\bremark\label{R3.1}
The above result can also be proved by the definition of Malliavin derivatives or approximation using the closeness of the operator $\mathbb{D}$.
\eremark

For any $v\in\R^d$, we aim to search for $h=h(v)\in {\rm Dom}\delta$ such that
\begin{equation}\label{3.2}
\nabla_v P(T)f(x)=\E(f(X^x(T))\delta(h)),\ f\in C_b^1(\R^d)
\end{equation}
holds.
We will show that $h$ satisfies \eqref{3.2} provided it is in ${\rm Dom}\delta$ with
\begin{equation}\label{3.3}
(R_H h)(t)=\int_0^t\sigma^{-1}(s)\left[\frac{T-s}{T}\nabla_v b(X^x(s))+\frac{1}{T}v\right]\d s, \ t\in[0,T].
\end{equation}

\btheorem\label{T3.1}
Assume (H1).
For $v\in\R^d$, let $h$ be given satisfying \eqref{3.3}.
If $h\in{\rm Dom}\delta$, then \eqref{3.2} holds.
\etheorem
\emph{Proof.}
By \eqref{3.1}, we know that the directional derivative process $\nabla_v X_\cdot^x$ satisfies the equation
$$\nabla_v X^x(t)=v+\int_0^t\nabla b(X^x(s))\nabla_v X^x(s)\d s, \ t\in[0,T].$$
On the other hand, observe that, for each $h\in\mathcal {H}, i=1,\cdots,d$,
\beqnn
\left\langle\sum_{j=1}^d(\sigma_{ij}\I_{[0,t]})e^j,h\right\rangle_\mathcal {H}
&=&\sum_{j=1}^d\langle K_H^*((\sigma_{ij}\I_{[0,t]})e^j),K_H^*h\rangle_{L^2([0,T],\R^d)}\\
&=&\sum_{j=1}^d\sum_{l=1}^d\int_0^T\left(K_H^*((\sigma_{ij}\I_{[0,t]})e^j)\right)^l(s)(K_H^*h)^l(s)\d s\\
&=&\sum_{j=1}^d\int_0^T\left(K_H^*(\sigma_{ij}\I_{[0,t]})\right)(s)(K_H^*h)^j(s)\d s\\
&=&\sum_{j=1}^d\int_0^T\int_s^T\sigma_{ij}(r)\I_{[0,t]}(r)\frac{\partial K_H}{\partial r}(r,s)(K_H^*h)^j(s)\d r\d s\\
&=&\sum_{j=1}^d\int_0^t\sigma_{ij}(r)\int_0^r\frac{\partial K_H}{\partial r}(r,s)(K_H^*h)^j(s)\d s\d r\\
&=&\sum_{j=1}^d\int_0^t\sigma_{ij}(r)\d (R_H h)^j(r),
\eeqnn
where the last equality follows from the fact: $R_H=K_H\circ K_H^\ast$.\\
Consequently, Lemma \ref{L3.1} leads to the fact that,
$$\langle\mathbb{D}X^{x,i}(t),h\rangle_\mathcal {H}=\sum_{j=1}^d\int_0^t(\nabla b(X^x(s)))_{ij}\langle\mathbb{D}X^{x,j}(s),h\rangle_\mathcal {H}\d s+\sum_{j=1}^d\int_0^t\sigma_{ij}(s)\d (R_H h)^j(s).$$
Set $g^i(t):=\nabla_v X^{x,i}(t)-\langle\mathbb{D}X^{x,i}(t),h\rangle_\mathcal {H}$, then it solves the equation
$$g^i(t)=v^i+\sum_{j=1}^d\int_0^t(\nabla b(X^x(s)))_{ij}g^j(s)\d s-\sum_{j=1}^d\int_0^t\sigma_{ij}(s)\d (R_H h)^j(s).$$
That is,
$$g(t)=v+\int_0^t\nabla b(X^x(s))g(s)\d s-\int_0^t\sigma(s)\d (R_H h)(s).$$
Therefore, if we set
$$(R_H h)(s)=\int_0^s\sigma^{-1}(r)\left[\nabla b(X^x(s))f(s)+\frac{1}{T}v\right]\d s,$$
then there holds $g(t)=\frac{T-t}{T}v$ and moreover $h$ satisfies \eqref{3.3}.
In particular, for each $i,g^i(T)=0,$ that is, $\nabla_v X^{x,i}(T)=\langle\mathbb{D}X^{x,i}(T),h\rangle_\mathcal {H}$.
As a consequence, for any $f\in C_b^1(\R^d)$, we obtain
\beqnn
\nabla_v P(T)f(x)&=&\E \nabla_vf(X^x(T))=\E((\nabla f)(X^x(T))\nabla_v X^x(T))
=\sum_{i=1}^d\E(\partial _if(X^x(T))\nabla_v X^{x,i}(T))\\
&=&\sum_{i=1}^d\E(\partial _if(X^x(T))\langle\mathbb{D}X^{x,i}(T),h\rangle_\mathcal {H})
=\E\langle\mathbb{D} f(X^x(T)),h\rangle_\mathcal {H}=\E(f(X^x(T))\delta(h)).
\eeqnn
\fin

To obtain explicit derivative formula from Theorem \ref{T3.1}, we need to calculate $\delta(h)$.
To this end, we need one more condition.

(H2) $\nabla b$ is H\"{o}lder continuous of order $1-1/(2H)<\beta_0\leq1$ with non-negative constant $L$:
\beqnn
|\nabla b(x)-\nabla b(y)|\leq L|x-y|^{\beta_0},\ \forall x,y\in\R^d.
\eeqnn

By fractional calculus, it is easy to get the following estimate.
\blemma\label{L3.2}
Assume (H1) and (H2).
Then, for any $s,t\in[0,T]$ there hold
\begin{eqnarray*}
|X_t-X_s|&\leq&C\Big\{\left[1+\e^{KT}(|x|+T)\right]|t-s|+\e^{KT}\|B^H\|_{\lambda_0}\left(T^{\lambda_0}+T^{\lambda_0+\alpha_0}\right)|t-s|\\
 &&~~~~+\|B^H\|_{\lambda_0}\left(|t-s|^{\lambda_0}+|t-s|^{\lambda_0+\alpha_0}\right)\Big\}\\
&\leq& C(T)(|t-s|+\|B^H\|_{\lambda_0}|t-s|^{\lambda_0}),
\end{eqnarray*}
where and in what follows, $C$ denotes a generic constant, $\lambda_0$ is chosen satisfying $1-\alpha_0<\lambda_0<H$ and $\lambda_0\beta_0>H-1/2$.
\elemma

\btheorem\label{T3.2}
Assume (H1) and (H2).
For $v\in\R^d$, let $h$ be given satisfying \eqref{3.3}.
Then there holds
\beqnn
\delta(h)
&=&\int_0^T\left\langle(K_H^*h)(t),\d W(t)\right\rangle\\
&=&\int_0^T\left\langle K_H^{-1}\left(\int_0^\cdot\sigma^{-1}(s)\left(\frac{T-s}{T}\nabla_v b(X^x(s))+\frac{1}{T}v\right)\d s\right)(t),\d W(t)\right\rangle.
\eeqnn
More precisely, $K_H^*h$ is shown explicitly in \eqref{T3.2add1} below.
\etheorem

\emph{Proof.}
According to \cite[Proposition 5.2.2]{Nualart06a}, we know that $h\in{\rm Dom}\delta$ if and only if $K_H^*h\in{\rm Dom}\delta_W$, where ${\rm Dom}\delta_W$ denotes the divergence operator with respect to the process $W$.
Next we want to show that $K_H^*h\in{\rm Dom}\delta_W$.
Due to  \cite[Proposition 1.3.11]{Nualart06a}, it suffices to prove $K_H^*h\in L^2_a([0,T]\times\Omega,\R^d)$, where  $L^2_a([0,T]\times\Omega,\R^d)$
represents the set of square integrable and adapted processes.
By \eqref{3.3}, we get
\begin{equation}\label{3.4}
(K_H^*h)(t)=K_H^{-1}\left(\int_0^\cdot\sigma^{-1}(s)\left(\frac{T-s}{T}\nabla_v b(X^x(s))+\frac{1}{T}v\right)\d s\right)(t), \ t\in[0,T].
\end{equation}
It is clear that the operator $K_H^{-1}$ preserves the adaptability property.
So, $(K_H^*h)$ is adapted.
Next, we focus on showing that $K_H^*h$ is square integrable.
By \eqref{2.1}, we obtain
\begin{eqnarray}\label{T3.2add1}
&&K_H^{-1}\left(\int_0^\cdot\sigma^{-1}(s)\left(\frac{T-s}{T}\nabla_v b(X^x(s))+\frac{1}{T}v\right)\d s\right)(t)\cr
&=&t^{H-\frac{1}{2}}D_{0+}^{H-\frac{1}{2}}\left[\cdot^{\frac{1}{2}-H}\sigma^{-1}(\cdot)\left(\frac{T-\cdot}{T}\nabla_v b(X^x(\cdot))+\frac{1}{T}v\right)\right](t)\cr
&=&\frac{1}{\Gamma(\frac{3}{2}-H)}\Bigg[t^{\frac{1}{2}-H}\sigma^{-1}(t)\left(\frac{T-t}{T}\nabla_v b(X^x(t))+\frac{1}{T}v\right)\cr
&&+\left(H-\frac{1}{2}\right)t^{H-\frac{1}{2}}\int_0^t\frac{t^{\frac{1}{2}-H}-r^{\frac{1}{2}-H}}{(t-r)^{\frac{1}{2}+H}}\sigma^{-1}(r)\left(\frac{T-r}{T}\nabla_v b(X^x(r))+\frac{1}{T}v\right)\d r\cr
&&+\left(H-\frac{1}{2}\right)\int_0^t\frac{\sigma^{-1}(t)-\sigma^{-1}(r)}{(t-r)^{\frac{1}{2}+H}}\left(\frac{T-t}{T}\nabla_v b(X^x(t))+\frac{1}{T}v\right)\cr
&&-\left(H-\frac{1}{2}\right)\frac{1}{T}\int_0^t(t-r)^{\frac{1}{2}-H}\sigma^{-1}(r)\nabla_v b(X^x(t))\d r\cr
&&+\left(H-\frac{1}{2}\right)\int_0^t\frac{T-r}{T}\sigma^{-1}(r)\frac{\nabla_v b(X^x(t))-\nabla_v b(X^x(r))}{(t-r)^{\frac{1}{2}+H}}\d r\Bigg]\cr
&=:&\frac{1}{\Gamma(\frac{3}{2}-H)}(I_1+I_2+I_3+I_4+I_5).
\end{eqnarray}
Note that
\begin{equation*}
\int_0^t\frac{t^{\frac{1}{2}-H}-r^{\frac{1}{2}-H}}{(t-r)^{\frac{1}{2}+H}}\d r
=\int_0^1\frac{\theta^{\frac{1}{2}-H}-1}{(1-\theta)^{\frac{1}{2}+H}}\d\theta\cdot t^{1-2H}=:C_0t^{1-2H},
\end{equation*}
where we make the change of variables $\theta=r/t,C_0$ is some positive constant.\\
This, together with (H1), leads to the fact that $I_i\in L^2_a([0,T]\times\Omega,\R^d), i=1,\cdots,4$.
As for $I_5$, by Lemma \ref{L3.2}, we deduce that
\beqnn
|I_5|\leq C\left(\int_0^t(t-r)^{\beta_0-1/2-H}\d r+\|B^H\|_{\lambda_0}^{\beta_0}\int_0^t(t-r)^{\lambda_0\beta_0-1/2-H}\d r\right).
\eeqnn
So, due to the Fernique theorem, $K_H^*h\in L^2([0,T]\times\Omega,\R^d)$ holds.\\
Again by \cite[Proposition 5.2.2]{Nualart06a}, we conclude that $\delta(h)=\delta_W(K_H^*h)=\int_0^T(K_H^*h)(t)\d W_t$.
This completes the proof.
\fin

As applications of the derivative formulae derived above, we will present Harnack type inequalities for $P(T)$.
To this end, we set, for all $s\in[0,T]$,
\beqnn
M_s=\int_0^s\left\langle K_H^{-1}\left(\int_0^\cdot\sigma^{-1}(r)\left(\frac{T-r}{T}\nabla_v b(X^x(r))+\frac{1}{T}v\right)\d r\right)(t),\d W(t)\right\rangle,
\eeqnn
and note that, by the Fernique theorem, there exists a positive constant $\theta_0$ such that $A_0:=\mathbb{E}\exp[\theta_0\|B^H\|_{\lambda_0}]<\infty$.

\btheorem\label{T3.3}
Assume (H1) and (H2).
Then there exist two positive constants $c'$ and $c''$ depending on $(T,K,\tilde{K},L,H)$
 such that for any $f\in\mathcal{B}_b^+(\R^d)$ and  $x,v\in\R^d,p>1$,
\beqnn
(P_T f(x))^p\leq P_Tf^p(x+v)\exp\left[\frac{p}{p-1}\left(c'+c''\left(1\vee\frac{p|v|}{p-1}\right)^{\frac{2\rho}{1-\rho}}\right)|v|^2\right]
\eeqnn
holds.
\etheorem

We first give an exponential estimate of $M$, which will be important in the proof of Theorem \ref{T3.3}.

\blemma\label{L3.3}
In the situation of Theorem \ref{T3.3}, for any $\lambda>0$, there exist constants $c'$ and $c''$ depending on $(T,K,\tilde{K},L,H)$ such that
\begin{equation*}
\E\exp\left[\frac{M_T}{\lambda}\right]
\leq\exp\left[\frac{|v|^2}{\lambda^2}\left(c'+c''\left(1\vee\frac{p|v|}{p-1}\right)^{\frac{2\beta_0}{1-\beta_0}}\right)\right].
\end{equation*}
\elemma

\emph{Proof.} Due to the fact that $\E\e^{L_t}\leq(\E\e^{2\langle L\rangle_t})^{1/2}$ for a continuous martingale $L_t$, we obtain
\begin{equation}\label{L3.3-1}
\E\e^{M_T/\lambda}\leq\left(\E\e^{2\langle M\rangle_T/\lambda^2}\right)^{1/2}.
\end{equation}
In view of the expression of the integrand of $M$ shown in \eqref{T3.2add1}, Lemma \ref{L3.2} and hypotheses imply that
\begin{eqnarray}\label{L3.3-2}
\langle M\rangle_T
&\leq&C\Bigg\{[(1+T)T^{-H}]^2+[(1+T)T^{\alpha_0-H}]^2+T^{2-2H}\cr
&&~~~~~+\int_0^T\left|\int_0^t\frac{|X^x(t)-X^x(r)|^{\beta_0}}{(t-r)^{\frac{1}{2}+H}}\d r\right|^2\d t\Bigg\}|v|^2\cr
&\leq&\left(a+b\|B^H\|_{\lambda_0}^{2\beta_0}\right)|v|^2,
\end{eqnarray}
where
\beqnn
a:=CT^{-2H}\left\{\left[1+\e^{KT}(|x|+T)\right]^2T^{2(\beta_0+1)}+(1+T)^2+[(1+T)T^{\alpha_0}]^2+T^2\right\}
\eeqnn
and
\beqnn
b:=CT^{2(1-H)}[\e^{2KT}(1+T^{\alpha_0})^2T^{2(\lambda_0+\beta_0)}+T^{2\lambda_0\beta_0}+T^{2(\lambda_0+\alpha_0)\beta_0}].
\eeqnn
Substituting \eqref{L3.3-2} into \eqref{L3.3-1} yields
\begin{equation}\label{L3.3-3}
\E\exp\left[\frac{M_T}{\lambda}\right]\leq\e^{(a/\lambda^2)|v|^2}\left(\E\exp\left[\frac{2b}{\lambda^2}\|B^H\|_{\lambda_0}^{2\beta_0}|v|^2\right]\right)^\frac{1}{2}.
\end{equation}
Note that
\beqnn
\frac{2b}{\lambda^2}\|B^H\|_{\lambda_0}^{2\beta_0}|v|^2
&=&\left[\frac{2b|v|^2}{\lambda^2}\left(\frac{|v|^2}{\lambda^2}\wedge1\right)^{-\beta_0}\left(\frac{\beta_0}{\theta_0}\right)^{\beta_0}\right]\cdot
\left[\left(\frac{|v|^2}{\lambda^2}\wedge1\right)^{\beta_0}\left(\frac{\theta_0}{\beta_0}\right)^{\beta_0}\|B^H\|_{\lambda_0}^{2\beta_0}\right]\cr
&\leq&(1-\beta_0)\left[2b\left(\frac{\beta_0}{\theta_0}\right)^{\beta_0}\right]^{\frac{1}{1-\beta_0}}\frac{1}
{\lambda^2}\left[|v|^2\vee\left(\frac{|v|}{\lambda^{\beta_0}}\right)^{\frac{2}{1-\beta_0}}\right]\cr
&&+\left(\frac{|v|^2}{\lambda^2}\wedge1\right)\theta_0\|B^H\|_{\lambda_0}^2\cr
&=:&\frac{2c}{\lambda^2}\left[|v|^2\vee\left(\frac{|v|}{\lambda^{\beta_0}}\right)^{\frac{2}{1-\beta_0}}\right]
+\left(\frac{|v|^2}{\lambda^2}\wedge1\right)\theta_0\|B^H\|_{\lambda_0}^2.
\eeqnn
Hence, by \eqref{L3.3-3} and the H\"{o}lder inequality we arrive at
\begin{eqnarray*}
\E\exp\left[\frac{M_T}{\lambda}\right]&\leq&\exp\left[\frac{a}{\lambda^2}|v|^2+\frac{c}{\lambda^2}\left(|v|^2\vee\left(\frac{|v|}{\lambda^{\beta_0}}\right)^{\frac{2}{1-\beta_0}}\right)\right]
\left(\E\exp\left[\left(\frac{|v|^2}{\lambda^2}\wedge1\right)\theta_0\|B^H\|_{\lambda_0}^2\right]\right)^{\frac{1}{2}}\cr
&\leq&
\exp\left[\frac{a}{\lambda^2}|v|^2+\frac{c}{\lambda^2}\left(|v|^2\vee\left(\frac{|v|}{\lambda^{\beta_0}}\right)^{\frac{2}{1-\beta_0}}\right)\right]
\left(\E\exp\left[\theta_0\|B^H\|_{\lambda_0}^2\right]\right)^{\frac{1}{2}\left(\frac{|v|^2}{\lambda^2}\wedge1\right)}\cr
&\leq&
\exp\left[\frac{|v|^2}{\lambda^2}\left(a+c\left(1\vee\frac{p|v|}{p-1}\right)^{\frac{2\beta_0}{1-\beta_0}}+\frac{\ln A_0}{2}\right)\right]\cr
&=:&\exp\left[\frac{|v|^2}{\lambda^2}\left(c'+c''\left(1\vee\frac{p|v|}{p-1}\right)^{\frac{2\beta_0}{1-\beta_0}}\right)\right].
\end{eqnarray*}
\fin

\emph{Proof of Theorem \ref{T3.3}.}
With the help of Theorem \ref{T3.1} and Theorem \ref{T3.2}, the argument of the proof is then standard.
We give its proof for the convenience of the reader.

Let $x,v\in\R^d$.
By Theorem \ref{T3.1} and the Young inequality \cite{Arnaudon&Thalmaier&Wang09a}, we derive that, for all $\lambda>0$,
\beqlb\label{3.5}
 |\nabla_v P(T)f(x)|\leq\lambda[P(T)(f\log f)(x)-(P(T)f)(x)(\log P(T)f)(x)]+\lambda P(T)f(x)\log\E\e^{\delta(h)/\lambda}.
\eeqlb
Now, let $\beta(s)=1+s(p-1), \gamma(s)=x+sv, s\in[0,1]$, we get, for any $p>1$,
\beqlb\label{3.6}
&&\frac{\d}{\d s}\log (P(T)f^{\beta(s)})^\frac{p}{\beta(s)}(\gamma(s))\cr
&=&\frac{p(p-1)}{\beta^2(s)}\frac{P(T)(f^{\beta(s)}\log f^{\beta(s)})-(P(T)f^{\beta(s)})\log P(T)f^{\beta(s)}}{P(T)f^{\beta(s)}}(\gamma(s))
+\frac{p}{\beta(s)}\frac{\nabla_v P(T)f^{\beta(s)}}{P(T)f^{\beta(s)}}(\gamma(s))\cr
&\geq&\frac{p}{\beta(s)P(T)f^{\beta(s)}(\gamma(s))}\bigg\{\frac{p-1}{\beta(s)}[P(T)(f^{\beta(s)}\log
f^{\beta(s)})(\gamma(s))\cr
&&-(P(T)f^{\beta(s)})\log P(T)f^{\beta(s)}(\gamma(s))]-|\nabla_v P(T)f^{\beta(s)}|(\gamma(s))\bigg\}\cr
&\geq&-\frac{p(p-1)}{\beta^2(s)}\log\E\e^{\delta(h)/\lambda},
\eeqlb
where in the last step we have applied \eqref{3.5} and then chosen $\lambda=\frac{p-1}{\beta(s)}$.\\
Applying Theorem \ref{T3.2} and Lemma \ref{L3.3} implies
\beqnn
\frac{\d}{\d s}\log (P(T)f^{\beta(s)})^\frac{p}{\beta(s)}(\gamma(s))
\geq-\frac{p}{p-1}\left[c'+c''\left(1\vee\frac{p|v|}{p-1}\right)^{\frac{2\beta_0}{1-\beta_0}}\right]|v|^2.
\eeqnn
Integrating on the interval $[0,1]$ with respect to $s$, we get the desired assertions.
\fin

\bcorollary\label{C3.1}
Assume (H1) and (H2).
Then
\beqnn
P(T)(\log f)(x)\leq \log P(T)f(x+v)+\left[c'+c''\left(1\vee\frac{p|v|}{p-1}\right)^{\frac{2\beta_0}{1-\beta_0}}\right]|v|^2
\eeqnn
holds for all $x,v\in\R^d$ and $f\in\mathcal{B}_b^+(\R^d)$.
\ecorollary
That is, log-Harnack inequalities hold.\\
In fact, since $\mathbb{R}^d$ is a length space, then by \cite[Proposition 2.2]{Wang10a} and Theorem \ref{T3.3}, we get the desired results.
In addition, combining \cite[Proposition 4.1]{DaPrato&Rockner&Wang09a} with Theorem \ref{T3.3}, we can deduce that $P(T)$ is strong Feller.

As for the case of $H<1/2$,
in \cite{Nualart&Ouknine02b} the authors proved the existence and uniqueness of the solution for the following stochastic equation with additive noise:
\begin{equation*}
\d X(t)=b(X(t))\d t+\d B^H(t), \ X_0=x\in\R^d,
\end{equation*}
Adopting the same arguments as $H>1/2$, we get the result for $<1/2$.
That is,

\btheorem\label{T3.4}
Assume (i) of (H1).
\begin{itemize}

\item[\rm(1)] For all $x,v\in\R^d$ and $f\in C_b^1(\R^d)$, there holds the Bismut formula
 \begin{eqnarray*}
&&\nabla_v P(T)f(x)=\E\left[f(X^x(T))\int_0^T\left\langle(K_H^*h)(t),\d W(t)\right\rangle\right]\\
&=&\E\left[f(X^x(T))\int_0^T\Bigg\langle\frac{t^{H-\frac{1}{2}}}{\Gamma(\frac{1}{2}-H)}\int_0^t
\frac{r^{\frac{1}{2}-H}}{(t-r)^{H+\frac{1}{2}}}\left(\frac{T-r}{T}\nabla_v b(X^x(r))+\frac{v}{T}\right)\d r,\d W(t)\Bigg\rangle\right].
\end{eqnarray*}

\item[\rm(2)] Harnack type inequalities
\beqnn
(P(T)f(x))^p\leq P(T)f^p(x+v)\exp\left[C\frac{p}{p-1}\frac{(1+T)^2}{T^{2H}}|v|^2\right]
\eeqnn
and
\beqnn
P(T)(\log f)(x)\leq \log P(T)f(x+v)+C\frac{(1+T)^2}{T^{2H}}|v|^2
\eeqnn
hold for all $x,v\in\R^d,p>1$ and $f\in\mathcal{B}_b^+(\R^d)$.

\end{itemize}
\etheorem

\section{Main results for SFDE}

\setcounter{equation}{0}

Let $r_0$ be fixed and denote $\mathscr{C}:=C([-r_0,0];\R^d)$ by the space of all continuous functions from
$[-r_0,0]$ to $\R^d$ equipped with the uniform norm $\|\phi\|_\infty=\sup_{-r_0\leq s\leq 0}|\phi(s)|$.
Let $T>r_0$ be a fixed number, for a map $g:[-r_0,T]\rightarrow\R^d$ and $t\geq0$, let $g_t\in\mathscr{C}$ be the
segment of $g(t)$. That is,
$$g_t(s)=g(t+s),\ s\in[-r_0,0].$$

In the present part, we study the following stochastic functional differential equation (SFDE)
\beqlb\label{4.1}
\d X(t)=b(X_t)\d t+\sigma(t)\d B^H(t), \ X_0=\xi\in\mathscr{C},
\eeqlb
where $b:\mathscr{C}\rightarrow\R^d,\sigma:[0,T]\rightarrow\R^d\times\R^d$ and $H>1/2$.

In \cite{Boufoussi&Hajji11}, the authors proved the existence and unqiueness result for the equation \eqref{4.1} under the Lipschitz conditions of $b$ and $\sigma$.
Let $X_t^\xi$ be the segment solution to \eqref{4.1} with $X_0=\xi$.
Define the Markov operator $(P_t)_{t\in[0,T]}$:
$$P_tf(\xi)=\E f(X_t^\xi),\ f\in\mathscr{B}_b(\mathscr{C}),\ \xi\in\mathscr{C}.$$
To establish derivative formula for $P_t$, we shall make use of the following assumption.
(A):
\begin{itemize}
\item[(i)] $b$ is Fr\'{e}chet differentiable such that $\nabla b$ is bounded and Lipschitz continuous;

\par

\item[(ii)] $\sigma^{-1}$ is H\"{o}lder continuous whose order is larger than $H-1/2$.
\end{itemize}
The main result of this part is the following.

\btheorem\label{T4.1}
Assume (A).
Then for any $C^1$-function $\gamma:[-r_0,T]\rightarrow\R$ satisfying $\gamma(t)=1$ for $t\in[-r_0,0]$ and $\gamma(t)=0$ for $t\in[T-r_0,T]$,
\begin{equation*}
\nabla_\eta P_tf(\xi)=\E\left[f(X_t^{\xi})\int_0^T\left\langle K_H^{-1}\left(\int_0^\cdot\sigma^{-1}(s)\left(\left(\nabla_{\Gamma_s}b\right)\left(X_s^\xi\right)-\gamma'(s)\eta(0)\right)\d s\right)(t),\d W(t)\right\rangle\right]
\end{equation*}
holds for all $\xi,\eta\in\mathscr{C}$ and $f\in C_b^1(\mathscr{C})$, where
$$\Gamma(t)=\left[\eta(t)\I_{[-r_0,0]}(t)+\eta(0)\I_{(0,T]}(t)\right]\gamma(t),\ t\in[-r_0,T].$$
\etheorem
The proof is modified from Section 3.
Let
$$\mathbb{D}_{R_Hh}X_t^\xi=\frac{\d}{\d\epsilon}\Big|_{\epsilon=0}X_t^\xi(w+\epsilon R_Hh),\ h\in\mathcal {H},$$
and
$$\nabla_\eta X_t^\xi=\frac{\d}{\d\epsilon}\Big|_{\epsilon=0}X_t^{\xi+\epsilon\eta},\ \eta\in\mathscr{C}.$$
Following the arguments of Lemma \ref{L3.1} and Theorem \ref{T3.1}, we conclude that
\begin{equation}\label{4.2}
\mathbb{D}_{R_Hh}X_t^\xi=\langle\mathbb{D}X_t^\xi,h\rangle_\mathcal{H}.
\end{equation}
Besides, in a standard way (see, for example, \cite[Lemma 3, Proposition 7]{Nualart&Saussereau09} or \cite[Theorem A.2]{Bao&Wang&Yuan13c} ) we derive the following result.
\blemma\label{L4.1}
Assume (A).
Then, $\mathbb{D}_{R_Hh}X^\xi(t)$ and $\nabla_\eta X^\xi(t)$ are the unique solutions to the equations
\begin{equation*}
\d Y(t)=\left[\left(\nabla_{Y_t}b\right)\left(X_t^\xi\right)+\sigma(t)(R_Hh)'(t)\right]\d t, \ Y_0=0,
\end{equation*}
 and
\begin{equation*}
\d Z(t)=\left(\nabla_{Z_t}b\right)\left(X_t^\xi\right)\d t, \ Z_0=\eta,
\end{equation*}
respectively.
\elemma

\emph{Proof of Theorem \ref{4.1}.}
By Lemma \ref{L4.1}, $V(t):=\nabla_\eta X^\xi(t)-\mathbb{D}_{R_Hh}X^\xi(t)$ solves the equation
\begin{equation*}
\d V(t)=\left[\left(\nabla_{V_t}b\right)\left(X_t^\xi\right)-\sigma(t)(R_Hh)'(t)\right]\d t, \ V_0=\eta.
\end{equation*}
Taking
$$R_Hh(t)=\int_0^t\sigma^{-1}(s)\left[\left(\nabla_{V_s}b\right)\left(X_s^\xi\right)-\gamma'(s)\eta(0)\right]\d s,\ t\in[0,T].$$
So, we get
$$V(t)=\left[\eta(t)\I_{[-r_0,0]}(t)+\eta(0)\I_{(0,T]}(t)\right]\gamma(t)=:\Gamma(t),\ t\in[-r_0,T].$$
In view of the definition of $\gamma$, we conclude that $V_T=0$, i.e., $\nabla_\eta X_T^\xi=\mathbb{D}_{R_Hh}X_T^\xi$.
Consequently,
$$R_Hh(t)=\int_0^t\sigma^{-1}(s)\left[\left(\nabla_{\Gamma_s}b\right)\left(X_s^\xi\right)-\gamma'(s)\eta(0)\right]\d s,\ t\in[0,T].$$
For $h$ satisfying the above relation, similar to the proof of Theorem \ref{T3.2} we may show that
$h\in{\rm Dom}\delta$ (or equivalently, $K_H^*h\in{\rm Dom}\delta_W$) holds.
Then for any $f\in C_b^1(\mathscr{C})$, applying \eqref{4.2} yields
\beqnn
\nabla_\eta P_T f(\xi)&=&\E\nabla_\eta f\left(X_T^\xi\right)=\E\left((\nabla f)\left(X_T^\xi\right)\nabla_\eta X_T^\xi\right)
=\E\left((\nabla f)\left(X_T^\xi\right)\mathbb{D}_{R_Hh}X_T^\xi\right)\\
&=&\E\left((\nabla f)\left(X_T^\xi\right)\langle\mathbb{D}X_T^\xi,h\rangle_\mathcal{H}\right)
=\E\left\langle\mathbb{D}f\left(X_T^\xi\right),h\right\rangle_\mathcal {H}=\E\left(f\left(X_T^\xi\right)\delta(h)\right)\\
&=&\E\left(f\left(X_T^\xi\right)
\int_0^T\left\langle K_H^{-1}\left(\int_0^\cdot\sigma^{-1}(s)\left(\left(\nabla_{\Gamma_s}b\right)\left(X_s^\xi\right)-\gamma'(s)\eta(0)\right)\d s\right)(t),\d W(t)\right\rangle\right).
\eeqnn
\fin

To conclude this section, we remark that, according to the proof of Theorem \ref{T3.3}, Harnack type inequalities for the operator $P_t$ associated with \eqref{4.1} also can be established..


\begin{thebibliography}{99}


\bibitem{Alos&Mazet&Nualart01a} E. Al\`{o}s, O. Mazet and D. Nualart,  {\it Stochastic calculus with respect to Gaussian processes}, Ann. Probab.
29(2001), 766--801.

\bibitem{Arnaudon&Thalmaier&Wang09a} M. Arnaudon, A. Thalmaier and F. Y. Wang, {\it Gradient estimates and Harnack inequalities on non-compact Riemannian
manifolds}, Stochastic Process. Appl. 119(2009), 3653-3670.



\bibitem{Bao&Wang&Yuan13c} J. Bao, F. Y. Wang and C. Yuan, {\it Bismut formulae and applications for functional SPDEs},
Bull. Sci. Math. 137(2013), 509--522.

\bibitem{Bao&Wang&Yuan13b} J. Bao, F. Y. Wang and C. Yuan,  {\it Derivative formula and Harnack inequality for degenerate functional SDEs},
Stoch. Dyn. 13(2013), 1--22.

\bibitem{Baudoin&Ouyang11} F. Baudoin and C. Ouyang, {\it Small-time kernel expansion for solutions of stochastic differential equations driven by fractional Brownian motions}, Stochastic Process Appl. 121(2011), 759--792.

\bibitem{Baudoin&Ouyang&Tindel11a} F. Baudoin, C. Ouyang and S. Tindel, {\it Upper bounds for the density of solutions of stochastic differential equations driven by fractional Brownian motions}, Ann. Inst. H. Poincar¨¦ Probab. Statist. 50(2014), 111--135.

\bibitem{Bismut84} J. M. Bismut, {\it Large Deviation and The Malliavin Calculus}, Birkh\"{a}user, MA, Boston, 1984.

\bibitem{Boufoussi&Hajji11} B. Boufoussi and S. Hajji,  {\it Functional differential equations driven by a fractional Brownian motion}, Comput. Math. Appl. 62(2011), 746--754.




\bibitem{Coutin&Qian02a} L. Coutin and Z. Qian, {\it Stochastic analysis, rough path analysis and fractional Brownian motions}, Probab. Theory Related Fields 122(2002), 108--140.



\bibitem{DaPrato&Rockner&Wang09a} G. Da Prato, M. R\"{o}ckner and F. Y. Wang, {\it Singular stochastic equations on Hilbert space: Harnack inequalities for their transition semigroups}, J. Funct. Anal. 257(2009), 992--1017.

\bibitem{Decreusefond&Ustunel98a} L. Decreusefond and A. S. \"{U}st\"{u}nel, {\it Stochastic analysis of the fractional Brownian motion}, Potential Anal. 10(1998), 177--214.

\bibitem{Diopa&Garrido-Atienza14} M. A. Diopa and M. J. Garrido-Atienza, {\it Retarded evolution systems driven by fractional Brownian
motion with Hurst parameter $H>1/2$}, Nonlinear Analysis 97(2014), 15--29.

\bibitem{Dong&Xie10} Z. Dong and Y. Xie,  {\it Ergodicity of linear SPDE Driven by L\'{e}vy noise},
J. Syst. Sci. Complex 23(2010), 137--152.




\bibitem{Elworthy&Li94} K. D. Elworthy and X. M. Li, {\it Formulae for the derivatives of heat semigroups}, J. Funct. Anal. 125(1994), 252--286.



\bibitem{Fan13} X. L. Fan, {\it Harnack inequality and derivative formula for SDE driven by fractional Brownian
motion}, Science in China-Mathematics 561(2013), 515--524.

\bibitem{Fan14} X. L. Fan, {\it Harnack-type inequalities and applications for SDE driven by fractional Brownian motion}, Stochastic Analysis and Applications 32(2014), 602--618.

\bibitem{Ferrante&Rovira06} M. Ferrante and C. Rovira,  {\it Stochastic delay differential equations driven by fractional Brownian motion with Hurst parameter $H>1/2$}, Bernoulli 12(2006), 85--100.



\bibitem{Guillin&Wang12} A. Guillin and F. Y. Wang,  {\it Degenerate Fokker-Planck equations: Bismut formula, gradient estimate and Harnack inequality},
J. Differential Equations 253(2012), 20--40.



\bibitem{Hairer05} M. Hairer, {\it Ergodicity of stochastic differential equations driven by fractional Brownian motion}, Ann. Probab. 33(2005), 703--758.

\bibitem{Hairer&Pillai11} M. Hairer and N. S. Pillai, {\it Ergodicity of hypoelliptic SDEs driven by fractional Brownian motion}, Ann. Inst. H. Poincar\'{e} Probab. Statist. 47(2011), 601--628.



\bibitem{Lyons98a} T. Lyons, {\it Differential equations driven by rough signals}, Rev. Mat. Iberoamericana 14(1998), 215--310.



\bibitem{Neuenkirch&Nourdin&Tindel08} A. Neuenkirch, I. Nourdin and S. Tindel,  {\it Delay equations driven by rough paths}, Electron. J. Probab. 13(2008), 2031--2068.

\bibitem{Nikiforov&Uvarov88} A. F. Nikiforov and V. B. Uvarov, {\it Special Functions of Mathematical Physics}, Birkh\"{a}user, Boston, 1988.

\bibitem{Nourdin&Simon06} I. Nourdin and T. Simon, {\it On the absolute continuity of one-dimensional SDEs driven by a fractional Brownian motion}, Statist. Probab. Lett. 76(2006), 907--912.

\bibitem{Nualart06a} D. Nualart, {\it The Malliavin Calculus and Related Topics, Second edition},
Springer-Verlag, Berlin, 2006.

\bibitem{Nualart&Ouknine02b} D. Nualart and Y. Ouknine,  {\it Regularization of differential equations by fractional noise},  Stochastic Process. Appl.
    102(2002), 103--116.

\bibitem{Nualart&Rascanu02a} D. Nualart and A. R\u{a}\c{s}canu, {\it Differential equations driven by fractional Brownian motion}, Collect. Math. 53(2002), 55--81.

\bibitem{Nualart&Saussereau09} D. Nualart and B. Saussereau, {\it Malliavin calculus for stochastic differential equations driven by a fractional Brownian motion}, Stochastic Process. Appl. 119(2009), 391--409.



\bibitem{Priola06} E. Priola,  {\it Formulae for the derivatives of degenerate diffusion semigroups},
J. Evol. Equ. 6(2006), 577--600.



\bibitem{Samko&Kilbas&Marichev} S. G. Samko, A. A. Kilbas and O. I. Marichev,  {\it Fractional Integrals and Derivatives, Theory and Applications}, Gordon and Breach Science Publishers, Yvendon, 1993.

\bibitem{Saussereau12} B. Saussereau, {\it Transportation inequalities for stochastic differential equations driven by a fractional Brownian motion}, Bernoulli 18(2012), 1--23.



\bibitem{Wang10a} F. Y. Wang, {\it Harnack inequalities on manifolds with boundary and applications}, J. Math. Pures Appl. 94(2010), 304--321.

\bibitem{Wang&Xu10a} F. Y. Wang and L. Xu, {\it Derivative formulae and applications for hyperdissipative stochastic Navier-Stokes/Burgers equations},
Infin. Dimens. Anal. Quantum Probab. Relat. Top. 15(2012), 1--19.

\bibitem{Wang&Zhang13a} F. Y. Wang and X. Zhang,  {\it Derivative formula and applications for degenerate diffusion semigroups},
J. Math. Pures Appl. 99(2013), 726--740.



\bibitem{Zahle98a} M. Z\"{a}hle, {\it Integration with respect to fractal functions and stochastic calculus I}, Probab. Theory Related Fields 111(1998), 333--374.

\bibitem{Zhang10a} X. Zhang, {\it Stochastic flows and Bismut formulas for stochastic Hamiltonian systems}, Stochastic Process. Appl. 120(2010), 1929--1949.

\end{thebibliography}
\end{document}